\documentclass[11pt,reqno,a4paper]{amsart}
\usepackage{amssymb}
\input amssym.def
\usepackage{amsmath,amsfonts,xcolor,mathtools,booktabs}
\usepackage[bookmarks=false]{hyperref}
\usepackage{amscd}
\usepackage[mathscr]{eucal}
\usepackage{enumitem}
\usepackage{orcidlink}
\allowdisplaybreaks
\numberwithin{equation}{section}

\hypersetup{colorlinks=true,citecolor={purple},linkcolor={teal},urlcolor={violet}}

\theoremstyle{plain}
\newtheorem{theorem}{Theorem}[section]
\newtheorem{corollary}[theorem]{Corollary}
\newtheorem{lemma}[theorem]{Lemma}

\theoremstyle{definition}

\newtheorem{remark}[theorem]{Remark}

\makeatletter
\@namedef{subjclassname@2020}{%
\textup{2020} Mathematics Subject Classification}
\makeatother



\begin{document}
\title[Gosper-type Lambert series identities of level $14$]{Gosper-type Lambert series identities of level $14$}	
\author[Russelle Guadalupe]{Russelle Guadalupe\orcidlink{0009-0001-8974-4502}}
\address{Institute of Mathematics, University of the Philippines Diliman\\
Quezon City 1101, Philippines}
\email{rguadalupe@math.upd.edu.ph}

\begin{abstract}
We derive two Gosper-type Lambert series identities of level $14$ which involve the $q$-constant $\Pi_q$ using a special case of Bailey's $_6\psi_6$ summation formula and certain propeties of $\eta$-quotients and generalized $\eta$-quotients on the congruence subgroup $\Gamma_0(14)$. 
\end{abstract}

\maketitle

\section{Introduction}\label{sec1}

Throughout this paper, we denote $q:=e^{2\pi i \tau}$ for an element $\tau$ in the complex upper half-plane $\mathbb{H}$. For any complex number $a$, we define $(a;q)_{\infty} := \prod_{n=0}^\infty (1-aq^n)$ and 
\begin{align*}
(a_1,a_2,\ldots,a_m;q)_\infty:= (a_1;q)_\infty (a_2;q)_\infty\cdots (a_m;q)_\infty. 
\end{align*}

In his seminal work on the convergence of power series, Lambert \cite{lamb} introduced an infinite series of the form 
\begin{align*}
\sum_{n\geq1}\dfrac{a_nq^n}{1-q^n}
\end{align*}
for a given sequence $\{a_n\}$ of complex numbers. This series is referred to as a Lambert series. There are a number of identities involving Lambert series, which include 
\begin{align*}
\dfrac{q}{(1-q)^2} &= \sum_{n\geq 1}\dfrac{\varphi(n)q^n}{1-q^n},\\
\sum_{n\geq 1}\sigma(n)q^n &=\sum_{n\geq 1}\dfrac{nq^n}{1-q^n} = \sum_{n\geq 1}\dfrac{q^n}{(1-q^n)^2},
\end{align*}
where $\varphi(n)$ is the number of positive integers at most $n$ that are relatively prime to $n$ and $\sigma(n)$ is the sum of all positive divisors of $n$.

In 2001, Gosper \cite{gosper} discovered without proof several identities involving the $q$-constant 
\begin{align*}
\Pi_q:=q^{1/4}\dfrac{(q^2;q^2)_\infty^2}{(q;q^2)_\infty^2},
\end{align*}
which is closely related to the $q$-analogues of the sine and cosine functions respectively given by 
\begin{align*}
\sin_q \pi z &:= q^{(z-1/2)^2}\prod_{n\geq 1}\dfrac{(1-q^{2n-2z})(1-q^{2n+2z-2})}{(1-q^{2n-1})^2},\\
\cos_q \pi z &:= q^{z^2}\prod_{n\geq 1}\dfrac{(1-q^{2n-2z-1})(1-q^{2n+2z+1})}{(1-q^{2n-1})^2}.
\end{align*}
One of the motivations for introducing $\Pi_q$ and the above $q$-analogues in his work is to derive some properties of the $q$-generalizations of hypergeometric identities. For instance, by exploiting some properties of $\sin_q \pi z$ and $\cos_q \pi z$, Gosper \cite{gosper} found a simple identity involving $\Pi_q$ given by
\begin{align*}
\dfrac{\Pi_q^2}{\Pi_{q^2}\Pi_{q^4}}-\dfrac{\Pi_{q^2}^2}{\Pi_{q^4}^2}=4,
\end{align*}
which can be used to devise a quadratically convergent method of calculating $\pi$. Gosper \cite{gosper} also found without proof the following Lambert series identities involving $\Pi_q$:
\begin{align}
&\left(\sum_{n\geq 1}\dfrac{q^n}{(1-q^n)^2}-2\sum_{n\geq 1}\dfrac{q^{2n}}{(1-q^{2n})^2}\right) = \dfrac{1}{24}\left(\dfrac{\Pi_q^4}{\Pi_{q^2}^2}-1\right)+\dfrac{2}{3}\Pi_{q^2}^2,\label{eq11}\\
&\dfrac{1}{\Pi_{q^3}^2}\left(\sum_{n\geq 1}\dfrac{q^{2n-1}}{(1-q^{2n-1})^2}-3\sum_{n\geq 1}\dfrac{q^{6n-3}}{(1-q^{6n-3})^2}\right) = \dfrac{\Pi_q}{\Pi_{q^3}},\label{eq12}\\
&\dfrac{1}{\Pi_{q^5}^2}\left(\sum_{n\geq 1}\dfrac{q^{2n-1}}{(1-q^{2n-1})^2}-5\sum_{n\geq 1}\dfrac{q^{10n-5}}{(1-q^{10n-5})^2}\right)= \sqrt{\dfrac{\Pi_q^3}{\Pi_{q^5}^3}-2\dfrac{\Pi_q^2}{\Pi_{q^5}^2}+5 \dfrac{\Pi_q}{\Pi_{q^5}}},\label{eq13}\\
&\phantom{\dfrac{1}{\Pi_{q^5}^2}\left(\sum_{n\geq 1}\dfrac{q^{2n-1}}{(1-q^{2n-1})^2}-5\sum_{n\geq 1}\dfrac{q^{10n-5}}{(1-q^{10n-5})^2}\right)}=\dfrac{\dfrac{\Pi_{q^5}}{\Pi_{q^{10}}}+16\dfrac{\Pi_{q^{10}}}{\Pi_{q^5}}}{\dfrac{\Pi_q}{\Pi_{q^5}}-4-\dfrac{\Pi_{q^5}}{\Pi_q}},\label{eq14}\\
&6\left(\sum_{n\geq 1}\dfrac{q^n}{(1-q^n)^2}-5\sum_{n\geq 1}\dfrac{q^{5n}}{(1-q^{5n})^2}\right)+1\nonumber\\
&\phantom{6\left(\sum_{n\geq 1}\dfrac{q^n}{(1-q^n)^2}\right)}=\left(\dfrac{\Pi_q}{\Pi_{q^5}}+2+5\dfrac{\Pi_{q^5}}{\Pi_q}\right)\left(\sum_{n\geq 1}\dfrac{q^{2n-1}}{(1-q^{2n-1})^2}-5\sum_{n\geq 1}\dfrac{q^{10n-5}}{(1-q^{10n-5})^2}\right),\label{eq15}\\
&\dfrac{1}{\Pi_{q^9}^2}\left(\sum_{n\geq 1}\dfrac{q^{2n-1}}{(1-q^{2n-1})^2}-9\sum_{n\geq 1}\dfrac{q^{18n-9}}{(1-q^{18n-9})^2}\right)\nonumber\\
&\phantom{\dfrac{1}{\Pi_{q^9}^2}\left(\sum_{n\geq1}\dfrac{q^{2n-1}}{(1-q^{2n-1})^2}-9\sum_{n\geq1}\right)}=\left(\dfrac{\Pi_q}{\Pi_{q^9}}+3\right)\sqrt{\left(\dfrac{\Pi_q}{\Pi_{q^9}}\right)^{3/2}-3\dfrac{\Pi_q}{\Pi_{q^9}}+3\left(\dfrac{\Pi_q}{\Pi_{q^9}}\right)^{1/2}},\label{eq16}\\
&3\left(\sum_{n\geq 1}\dfrac{q^n}{(1-q^n)^2}-9\sum_{n\geq 1}\dfrac{q^{9n}}{(1-q^{9n})^2}\right)+1\nonumber\\
&\phantom{3\left(\sum_{n\geq 1}\dfrac{q^n}{(1-q^n)^2}\right)}=\left(\sqrt{\dfrac{\Pi_q}{\Pi_{q^9}}}+3\sqrt{\dfrac{\Pi_{q^9}}{\Pi_q}}\right)\left(\sum_{n\geq 1}\dfrac{q^{2n-1}}{(1-q^{2n-1})^2}-9\sum_{n\geq 1}\dfrac{q^{18n-9}}{(1-q^{18n-9})^2}\right).\label{eq17}	
\end{align}

El Bachraoui \cite{bach} obtained (\ref{eq11}) and (\ref{eq12}) using Gosper's $q$-trigonometry. He \cite{he1} proved (\ref{eq13}), (\ref{eq14}), and (\ref{eq15}) using Ramanujan's modular equations of degree five. He \cite{he2} also proved (\ref{eq16}) using the theory of modular forms. Wang \cite{wang} established the above identities by observing that these can be interpreted as equalities of certain modular forms on the congruence subgroup $\Gamma_0(N)$, and called such equality an identity of level $N$. Thus, (\ref{eq11}) is an identity of level $4$, (\ref{eq12}) is an identity of level $6$, (\ref{eq13}) and (\ref{eq15}) are identities of level $10$, (\ref{eq16}) and (\ref{eq17}) are identities of level $18$, and (\ref{eq14}) is an identity of level $20$. Wang \cite{wang} also found another Gosper-type Lambert series identities of level $18$. Yathirajsharma, Harshitha, and Vasuki \cite{yhv} confirmed the above identities using Ramanujan's theta function identities and proved two new Lambert series identities of level $20$. Recently, Yathirajsharma \cite{yat} provided new proofs of (\ref{eq14}) and an equivalent form of (\ref{eq16}) by employing some $q$-trigonometric identities. 

We continue in this paper our investigation on Gosper-type Lambert series identities by establishing two new such identities of level $14$. We first set 
\begin{align*}
z &:= \dfrac{1}{\Pi_{q^7}^2}\left(\sum_{n\geq 1}\dfrac{q^{2n-1}}{(1-q^{2n-1})^2}-7\sum_{n\geq 1}\dfrac{q^{14n-7}}{(1-q^{14n-7})^2}\right),\\
w &:= \dfrac{1}{\Pi_{q^7}^2}\left[4\left(\sum_{n\geq 1}\dfrac{q^n}{(1-q^n)^2}-7\sum_{n\geq 1}\dfrac{q^{7n}}{(1-q^{7n})^2}\right)+1\right].
\end{align*}

\begin{theorem}\label{thm11} We have the identity
\begin{align}\label{eq18}
z^3+4\dfrac{\Pi_q}{\Pi_{q^7}} z^2-3\dfrac{\Pi_q^2}{\Pi_{q^7}^2} z-\dfrac{\Pi_q}{\Pi_{q^7}}\left(\dfrac{\Pi_q^4}{\Pi_{q^7}^4}+4\dfrac{\Pi_q^2}{\Pi_{q^7}^2}+49\right)=0.
\end{align}
\end{theorem}

\begin{theorem}\label{thm12} Let $f:= w/z$. Then $f$ satisfies the identity
\begin{align}\label{eq19}
&\dfrac{\Pi_q^2}{\Pi_{q^7}^2}\left(\dfrac{\Pi_q^4}{\Pi_{q^7}^4}+4\dfrac{\Pi_q^2}{\Pi_{q^7}^2}+49\right)f^3 - \dfrac{\Pi_q^2}{\Pi_{q^7}^2}\left(2\dfrac{\Pi_q^4}{\Pi_{q^7}^4}+5\dfrac{\Pi_q^2}{\Pi_{q^7}^2}+98\right)f^2\nonumber\\
&-2\dfrac{\Pi_q^2}{\Pi_{q^7}^2}\left(5\dfrac{\Pi_q^4}{\Pi_{q^7}^4}+22\dfrac{\Pi_q^2}{\Pi_{q^7}^2}+245\right)f-\left(\dfrac{\Pi_q^2}{\Pi_{q^7}^2}-4\dfrac{\Pi_q}{\Pi_{q^7}}+7\right)^2\left(\dfrac{\Pi_q^2}{\Pi_{q^7}^2}+4\dfrac{\Pi_q}{\Pi_{q^7}}+7\right)^2=0.
\end{align}
\end{theorem}

We observe that the level $14$ identities (\ref{eq18}) and (\ref{eq19}) are not as simple as that of identities (\ref{eq11}) -- (\ref{eq17}). This is simply due to the fact that the underlying congruence subgroup $\Gamma_0(14)$ associated with the level $14$ identities has genus one (as opposed to $\Gamma_0(N)$ which has genus zero for $N\in\{4,6,10,18\}$), which will be essential in establishing these identities.

We organize the rest of the paper as follows. In Section \ref{sec2}, we recall important identities, including a special case of Bailey's well-poised $_6\psi_6$ summation formula \cite[(4.7)]{bailey1}, which will be needed to prove Theorems \ref{thm11} and \ref{thm12}. We also recall basic facts on modular functions, particularly $\eta$-quotients and generalized $\eta$-quotients, on congruence subgroups $\Gamma_1(N)$ and $\Gamma_0(N)$. In Section \ref{sec3}, we demonstrate Theorems \ref{thm11} and \ref{thm12}, respectively, by employing Bailey's formula and generalized $\eta$-quotients. We have performed most of our computations via \textit{Mathematica}. We remark that our methods are different from that of He \cite{he1,he2} and Wang \cite{wang} in that we only involve modular functions instead of using Eisenstein series and constructing auxiliary modular forms of positive weight. We also remark that our methods, in principle, can be applied to any congruence subgroup $\Gamma_0(N)$ that will give rise to a Lambert series identity of level $N$.

\section{Preliminaries}\label{sec2}

We first recall the Ramanujan's general theta function 
\begin{align*}
f(a,b) =\sum_{n=-\infty}^\infty a^{n(n+1)/2}b^{n(n-1)/2}=(-a,-b,ab;ab)_\infty
\end{align*}
which holds for any complex numbers $a$ and $b$ with $|ab| < 1$. The last equality follows from the Jacobi triple product identity \cite[p. 35, Entry 19]{berndt3}. We define $f(-q):=f(-q,-q^2)=(q;q)_\infty$. We require the following identity due to Bailey, which relates Lambert series with products of theta functions.  

\begin{lemma}[\cite{bailey2}]\label{lem21}
For $|ab| < 1$, we have
\begin{align*}
\sum_{n=-\infty}^\infty\left[\dfrac{aq^n}{(1-aq^n)^2}-\dfrac{bq^n}{(1-bq^n)^2}\right]=af^6(-q)\dfrac{f(-ab,-\frac{q}{ab})f(-\frac{b}{a},-\frac{aq}{b})}{f^2(-a,-\frac{q}{a})f^2(-b,-\frac{q}{b})}.
\end{align*}
\end{lemma}

\begin{remark}\label{rem22}
Bailey \cite{bailey2} proved Lemma \ref{lem21} using the theory of elliptic functions. It can be shown that Lemma \ref{lem21} follows from his very well-poised $_6\psi_6$ summation formula \cite[(4.7)]{bailey1}.
\end{remark}

We next recall basic facts about modular functions on congruence subgroups
\begin{align*}
\Gamma_1(N) &:= \left\lbrace\begin{bmatrix}
a & b\\ c& d
\end{bmatrix} \in \mathrm{SL}_2(\mathbb{Z}) : a-1\equiv c\equiv d-1\equiv 0\pmod N\right\rbrace,\\
\Gamma_0(N) &:= \left\lbrace\begin{bmatrix}
a & b\\ c& d
\end{bmatrix} \in \mathrm{SL}_2(\mathbb{Z}) : c\equiv 0\pmod N\right\rbrace.
\end{align*}

Recall that an element on a congruence subgroup $\Gamma$ acts on the extended upper half-plane $\mathbb{H}^\ast := \mathbb{H}\cup \mathbb{Q}\cup \{\infty\}$ by a linear fractional transformation. Under this action, we call the equivalence classes of $\mathbb{Q}\cup \{\infty\}$ cusps, and there are finitely many inequivalent cusps of $\Gamma$ by \cite[Proposition 1.32]{shim}. We call a meromorphic function $f:\mathbb{H}\rightarrow\mathbb{C}$ a modular function on $\Gamma$ if 
\begin{enumerate}
\item[(1)] for all $\gamma\in \Gamma$, $f(\gamma\tau)=f(\tau)$, and 
\item[(2)] for every cusp $r$ of $\Gamma$ and $\gamma\in\mathrm{SL}_2(\mathbb{Z})$ with $\gamma(\infty)=r$, we have the $q$-expansion
\begin{align*}
f(\gamma\tau) = \sum_{n\geq n_0} a_nq^{n/h}
\end{align*}
for some integers $h$ and $n_0$ with $a_{n_0}\neq 0$. Here, $n_0$ is the order of vanishing of $f(\tau)$ at $r$, denoted by $\mathrm{ord}_r\,f(\tau)$, and $h$ is the width of $r$, which is the smallest positive integer such that $\gamma[\begin{smallmatrix}
1 & h\\0 & 1
\end{smallmatrix}]\gamma^{-1}\in \pm\Gamma$. We say that $f(\tau)$ has a zero (resp., a pole) at $r$ if $\mathrm{ord}_r\,f(\tau) > 0$ (resp.,  $\mathrm{ord}_r\,f(\tau) < 0$).
\end{enumerate}

The following results describe the inequivalent cusps of $\Gamma_0(N)$ and $\Gamma_1(N)$ with their respective widths.

\begin{lemma}\label{lem23}
Let $a,a',c$ and $c'$ be integers with $\gcd(a,c)=\gcd(a',c')=1$. Let $S_{\Gamma_0(N)}$ be the set of all inequivalent cusps of $\Gamma_0(N)$. We denote $\pm1/0$ as $\infty$. Then 
\begin{enumerate}
\item[(1)] $a/c$ and $a'/c'$ are equivalent over $\Gamma_0(N)$ if and only if there exist an integer $n$ and an element $\overline{s}\in (\mathbb{Z}/N\mathbb{Z})^\times$ such that $(a',c')\equiv(\overline{s}^{-1}a+nc,\overline{s}c)\pmod N$, and 
\item[(2)] we have 
\begin{align*}
	S_{\Gamma_0(N)} = \{a_{c,j}/c\in \mathbb{Q} : 0 < c \mid N, 0 < a_{c,j}\leq N, \gcd(a_{c,j},N)=1,\\ a_{c,j}= a_{c,j'}\stackrel{\text{def}}{\iff} a_{c,j}\equiv a_{c,j'}\pmod{\gcd(c,N/c)}\},
\end{align*}
and the width of the cusp $a/c$ in $\Gamma_0(N)$ is $N/\gcd(c^2,N)$.
\end{enumerate}
\end{lemma}

\begin{proof}
See \cite[Corollary 4(1)]{chokoop}.
\end{proof}

\begin{lemma}\label{lem24}
Let $a,a',c$ and $c'$ be integers with $\gcd(a,c)=\gcd(a',c')=1$. Let $S_{\Gamma_1(N)}$ be the set of all inequivalent cusps of $\Gamma_1(N)$. We denote $\pm1/0$ as $\infty$. Then 
\begin{enumerate}
\item[(1)] $a/c$ and $a'/c'$ are equivalent under $\Gamma_1(N)$ if and only if there exists $n\in\mathbb{Z}$ such that $[\begin{smallmatrix}
	a' \\ c'
\end{smallmatrix}]\equiv \pm[\begin{smallmatrix}
	a+nc \\ c
\end{smallmatrix}]\pmod{N}$,
\item[(2)] we have
\begin{align*}
	S_{\Gamma_1(N)} =&\{y_{c,j}/x_{c,i}\in\mathbb{Q} : 0 < c\mid N, 0<s_{c,i},a_{c,j}\leq N,\\ 
	&\gcd(s_{c,i},N)=\gcd(a_{c,j},N)=1,\\
	& s_{c,i}=s_{c,i'}\stackrel{\text{def}}{\iff} s_{c,i}\equiv \pm s_{c,i'}\pmod{N/c},\\
	& a_{c,j}=a_{c,j'}\stackrel{\text{def}}{\iff} a_{c,j}\equiv \begin{cases}
		\pm a_{c,j'}\pmod{c} &\text{ if } c\in\{N/2,N\},\\
		a_{c,j'}\pmod{c} &\text{ otherwise, }
	\end{cases}\\
	&(x_{c,i},y_{c,j})\in\mathbb{Z}^2, \gcd(x_{c,i},y_{c,j})=1,\\
	&x_{c,i}\equiv cs_{c,i}\pmod{N}, y_{c,j}\equiv a_{c,j}\pmod{N}\},
\end{align*}
and the width of the cusp $a/c$ in $\Gamma_1(N)$ is $1$ (respectively, $N/\gcd(c,N)$) if $N=4$ and $\gcd(c,4)=2$ (respectively, otherwise).
\end{enumerate}
\end{lemma}

\begin{proof}
See \cite[Corollary 4(2)]{chokoop}.
\end{proof}

We denote $\mathcal{M}^\infty(N)$ as the set of all modular functions on $\Gamma_0(N)$ that have poles only at infinity. The next result concerns about a polynomial relation between $x, y\in \mathcal{M}^\infty(N)$ whose orders of vanishing at infinity are coprime, provided that $\Gamma_0(N)$ has a nonzero genus.

\begin{lemma}\label{lem25}
Suppose $\Gamma_0(N)$ has a nonzero genus. Let $x, y\in \mathcal{M}^\infty(N)$ with unique pole of orders of vanishing $m$ and $n$, respectively, at infinity with $\gcd(m,n)=1$, such that the leading coefficients in their $q$-expansions are both $1$. Then $x$ and $y$ generate $\mathcal{M}^\infty(N)$, and there is a polynomial 
\begin{align*}
F(X,Y) = X^n-Y^m +\sum_{\substack{a, b\geq 0\\ am+bn\leq mn}}C_{a,b}X^aY^b\in \mathbb{C}[X,Y]
\end{align*}
such that $F(x,y)=0$. In addition, $F(x,Y)$ is the minimal polynomial for $y$ over $\mathbb{C}[x]$.
\end{lemma}

\begin{proof}
See \cite[Lemma 1, Lemma 2]{yang2}.
\end{proof}

We next define an $\eta$-quotient, which is a function of the form
\begin{align*}
f(\tau) = \prod_{\delta\mid N} \eta^{r_{\delta}}(\delta\tau)
\end{align*}
for some indexed set $\{r_\delta\in\mathbb{Z} : \delta\mid N\}$, where $\eta(\tau)=q^{1/24}(q;q)_\infty$ is the Dedekind eta function. The following results give conditions for an $\eta$-quotient to be a modular function on $\Gamma_0(N)$. 

\begin{lemma}\label{lem26}
Let $f(\tau) = \prod_{\delta\mid N} \eta^{r_{\delta}}(\delta\tau)$ be an $\eta$-quotient with $\sum_{\delta\mid N} r_{\delta}=0$. Then $f$ is a modular function on $\Gamma_0(N)$ with character
\begin{align*}
\chi(d):=\left(\dfrac{\prod_{\delta\mid N} \delta^{r_\delta}}{d}\right),
\end{align*}
where $(\frac{\cdot}{d})$ is the Kronecker symbol, if and only if
\begin{equation*}
\sum_{\delta\mid N} \delta r_{\delta} \equiv 0\pmod{24}\quad\text{ and }\quad\sum_{\delta\mid N} \dfrac{N}{\delta}r_{\delta} \equiv 0\pmod {24}.
\end{equation*}
\end{lemma} 

\begin{proof}
See \cite[Proposition 5.9.2]{cohen}.
\end{proof}

\begin{lemma}\label{lem27}
Let $c, d$, and $N$ be positive integers with $d\mid N$ and $\gcd(c,d)=1$ and let $f(\tau) = \prod_{\delta\mid N} \eta^{r_{\delta}}(\delta\tau)$ be an $\eta$-quotient satisfying the conditions of Lemma \ref{lem26}. Then the order of vanishing of $f(\tau)$ at the cusp $c/d$ of $\Gamma_0(N)$ is 
\begin{align*}
\mathrm{ord}_{c/d}\, f(\tau) = \dfrac{N}{24d\gcd(d,\frac{N}{d})}\sum_{\delta\mid N} \gcd(d,\delta)^2\cdot\dfrac{r_{\delta}}{\delta}.
\end{align*}
\end{lemma}

\begin{proof}
See \cite[Lemma 4.1.3]{aygin}.
\end{proof}

We introduce generalized $\eta$-quotients due to Yang \cite{yang1}, which are functions of the form
\begin{align*}
f(\tau) = \prod_{1\leq g\leq \lfloor N/2\rfloor}\eta_{N,g}^{r_g}(\tau)
\end{align*}
for some indexed set $\{r_g \in\mathbb{Z}: 1\leq g\leq \lfloor N/2\rfloor\}$, where
\begin{align*}
\eta_{N, g}(\tau)=q^{NB_2(g/N)/2}(q^g,q^{N-g};q^N)_\infty
\end{align*} 
is the generalized Dedekind eta function with $B_2(t) := t^2-t+1/6$. The following results give the transformation formula of $\eta_{N, g}$ on $\Gamma_0(N)$ and conditions for a generalized $\eta$-quotient to be modular on $\Gamma_1(N)$.

\begin{lemma}\label{lem28}
The function $\eta_{N,g}(\tau)$ satisfies $\eta_{N,g+N}(\tau)=\eta_{N,-g}(\tau)=-\eta_{N,g}(\tau)$. Moreover, let $\gamma = [\begin{smallmatrix}
a & b\\ cN & d
\end{smallmatrix}]\in \Gamma_0(N)$ with $c\neq 0$. Then we have 
\begin{align*}
\eta_{N,g}(\gamma\tau) = \varepsilon(a,bN,c,d)e^{\pi i(g^2ab/N - gb)}\eta_{N,ag}(\tau),
\end{align*}
where 
\begin{align*}
\varepsilon(a,b,c,d) = \begin{cases}
	e^{\pi i(bd(1-c^2)+c(a+d-3))/6}, & \text{ if }c\equiv 1\pmod{2},\\
	-ie^{\pi i(ac(1-d^2)+d(b-c+3))/6}, & \text{ if }c\equiv 0\pmod{2}.
\end{cases}
\end{align*}
\end{lemma}

\begin{proof}
See \cite[Corollary 2]{yang1}.
\end{proof}

\begin{lemma}\label{lem29}
Suppose $f(\tau)=\prod_{1\leq g\leq \lfloor N/2\rfloor}\eta_{N,g}^{r_g}(\tau)$ is a generalized $\eta$-quotient such that
\begin{enumerate}
\item[(1)] $\sum_{1\leq g\leq \lfloor N/2\rfloor}r_g\equiv 0\pmod{12}$, 
\item[(2)] $\sum_{1\leq g\leq \lfloor N/2\rfloor}gr_g\equiv 0\pmod{2}$, and 
\item[(3)] $\sum_{1\leq g\leq \lfloor N/2\rfloor}g^2r_g\equiv 0\pmod{2N}$.
\end{enumerate}
Then $f(\tau)$ is a modular function on $\Gamma_1(N)$. 
\end{lemma}

\begin{proof}
See \cite[Corollary 3]{yang1}.
\end{proof}

\begin{lemma}\label{lem210}
Let $N$ be a positive integer and $\gamma=[\begin{smallmatrix}
a & b\\ c& d
\end{smallmatrix}]\in\mathrm{SL}_2(\mathbb{Z})$. Then the first term of the $q$-expansion of $\eta_{N, g}(\gamma\tau)$ is $\varepsilon q^{\delta}$, where $|\varepsilon|=1$ and 
\[\delta = \dfrac{\gcd(c,N)^2}{2N}P_2\left(\dfrac{ag}{\gcd(c,N)}\right)\]
with $P_2(t):=B_2(\{t\})$ and $\{t\}$ is the fractional part of $t$.
\end{lemma}

\begin{proof}
See \cite[Lemma 2]{yang1}.
\end{proof}

\begin{corollary}\label{cor211}
If $f(\tau)=\prod_{1\leq g\leq \lfloor N/2\rfloor}\eta_{N,g}^{r_g}(\tau)$ is a generalized $\eta$-quotient on $\Gamma_1(N)$ with $N > 4$, then the order of vanishing of $f(\tau)$ at the cusp $a/c$ of $\Gamma_1(N)$ is 
\begin{align*}
\mathrm{ord}_{a/c}\, f(\tau) = \dfrac{\gcd(c,N)}{2}\sum_{1\leq g\leq \lfloor N/2\rfloor} P_2\left(\dfrac{ag}{\gcd(c,N)}\right)r_g.
\end{align*}
\end{corollary}

\begin{proof}
This follows from Lemmas \ref{lem24} and \ref{lem210}.
\end{proof}

\section{Proofs of Theorems \ref{thm11} and \ref{thm12}}\label{sec3}

We prove in this section Theorems \ref{thm11} and \ref{thm12} by considering certain sums of generalized $\eta$-quotients on the genus one congruence subgroup $\Gamma_0(14)$. 

\begin{proof}[Proof of Theorem \ref{thm11}]
Using Lemmas \ref{lem26} and \ref{lem27}, we find that the $\eta$-quotient
\begin{align*}
g^2 = \dfrac{\eta^8(2\tau)\eta^4(7\tau)}{\eta^4(\tau)\eta^8(14\tau)}
\end{align*}	
where $g:=g(\tau)=\Pi_q/\Pi_{q^7}$, is in $\mathcal{M}^\infty(14)$ with $\mathrm{ord}_\infty\,g^2=-3$ and $\mathrm{ord}_{1/2}\,g^2=3$. From Lemma \ref{lem21}, we replace $q$ with $q^{14}$ and set $b = q^7$ to get
\begin{align}
\sum_{n=-\infty}^\infty\left[\dfrac{aq^{14n}}{(1-aq^{14n})^2}-\dfrac{q^{14n+7}}{(1-q^{14n+7})^2}\right]=aq^{-7/2}\Pi_{q^7}^2\dfrac{f^2(-aq^7,-\frac{q^7}{a})}{f^2(-a,-\frac{q^{14}}{a})}. \label{eq31}
\end{align}
Setting $a=q, q^3, q^5$ in (\ref{eq31}) and adding the resulting identities, we have 
\begin{align*}
z = g_1(\tau)+g_2(\tau)+g_3(\tau),
\end{align*}
where 
\begin{align*}
g_1(\tau) := \dfrac{\eta_{14, 6}^2(\tau)}{\eta_{14, 1}^2(\tau)},\qquad	g_2(\tau) := \dfrac{\eta_{14, 4}^2(\tau)}{\eta_{14, 3}^2(\tau)},\qquad g_3(\tau) := \dfrac{\eta_{14, 2}^2(\tau)}{\eta_{14, 5}^2(\tau)}.
\end{align*}
Note that $g_1(\tau)g_2(\tau)g_3(\tau) = g$, so that
\begin{align*}
gz = g_1^2(\tau)g_2(\tau)g_3(\tau) + g_1(\tau)g_2^2(\tau)g_3(\tau) + g_1(\tau)g_2(\tau)g_3^2(\tau).
\end{align*}
We see from Lemma \ref{lem29} that $g_1(\tau)g_2(\tau), g_2(\tau)g_3(\tau)$, and $g_1(\tau)g_3(\tau)$ are modular on $\Gamma_1(14)$. Furthermore, Lemma \ref{lem28} shows that
\begin{align*}
	g_1(\gamma\tau) = -g_2(\tau), \qquad g_2(\gamma\tau)=-g_3(\tau),\qquad g_3(\gamma\tau) = -g_1(\tau)
\end{align*}
where $\gamma:=[\begin{smallmatrix}
	3 & 1\\ 14 & 5
\end{smallmatrix}]$. As $\Gamma_0(14)$ is generated by $\Gamma_1(14)$ and $\gamma$, we infer that the functions 
\begin{align*}
	f_1 &:= g_1(\tau)g_2(\tau)+g_1(\tau)g_3(\tau)+g_2(\tau)g_3(\tau)
\end{align*}
and $gz$ are modular on $\Gamma_0(14)$. We next apply Corollary \ref{cor211} to find the orders of vanishing of $f_1$ and $gz$ at each element of the set
\begin{align*}
	S_{\Gamma_1(14)} = \{\infty, 0,1/3,1/5,1/2,1/6,1/10,1/7,3/7,5/7,3/14,5/14\}
\end{align*}
of inequivalent cusps of $\Gamma_1(14)$ obtained from Lemma \ref{lem24}. Observe that for $r\in S_{\Gamma_1(14)}$, we have
\begin{align*}
	\mathrm{ord}_r\,f_1 &\geq \min\{\mathrm{ord}_r\, g_1(\tau)g_2(\tau),\mathrm{ord}_r\, g_1(\tau)g_3(\tau), \mathrm{ord}_r\, g_2(\tau)g_3(\tau)\},\\
	\mathrm{ord}_r\,gz &\geq \min\{\mathrm{ord}_r\, g_1^2(\tau)g_2(\tau)g_3(\tau),\mathrm{ord}_r\, g_1(\tau)g_2^2(\tau)g_3(\tau), \mathrm{ord}_r\,g_1(\tau)g_2(\tau)g_3^2(\tau)\}.
\end{align*}
Equality holds in the first inequality when $\mathrm{ord}_r\, g_i(\tau)g_j(\tau)$ are not all equal for $1\leq i < j \leq 3$ and in the second inequality when $\mathrm{ord}_r\, g_1^2(\tau)g_2(\tau)g_3(\tau), \mathrm{ord}_r\, g_1(\tau)g_2^2(\tau)g_3(\tau)$, and $\mathrm{ord}_r g_1(\tau)g_2(\tau)g_3^2(\tau)$ are not all equal. We tabulate the values of these orders as shown in Table \ref{tab:tbl31}.

\begin{table}[h]
	\caption{The orders of vanishing of $f_1$ and $gz$ at the cusps $r$ of $\Gamma_1(14)$}\label{tab:tbl31}%
	\begin{tabular}{@{}cccccc@{}}
		\toprule
		cusp $r$ & $\infty$ & $3/14$ & $5/14$ & $1/2,1/6,1/10$ & any other $r$  \\
		\midrule
		$\mathrm{ord}_r\, g_1(\tau)g_2(\tau)$ & $-3$ & $1$ & $-1$ & $1$ & $0$\\
		\midrule
		$\mathrm{ord}_r\, g_1(\tau)g_3(\tau)$ & $-1$ & $-3$ & $1$& $1$ & $0$\\
		\midrule
		$\mathrm{ord}_r\, g_2(\tau)g_3(\tau)$ & $1$ & $-1$ & $-3$ & $1$ & $0$ \\
		\midrule
		$\mathrm{ord}_r\, g_1^2(\tau)g_2(\tau)g_3(\tau)$ & $-4$ & $-2$ & $0$& $2$ & $0$\\
		\midrule
		$\mathrm{ord}_r\, g_1(\tau)g_2^2(\tau)g_3(\tau)$ & $-2$ & $0$ & $-4$ & $2$ & $0$\\
		\midrule
		$\mathrm{ord}_r\, g_1(\tau)g_2(\tau)g_3^2(\tau)$ & $0$ & $-4$ & $-2$ & $2$ & $0$\\
		\midrule
		$\mathrm{ord}_r\, f_1$ & $-3$ & $-3$ & $-3$ & $\geq 1$ & $\geq 0$\\
		\midrule
		$\mathrm{ord}_r\, gz$ & $-4$ & $-4$ & $-4$ & $\geq 2$ & $\geq 0$\\
		\bottomrule
	\end{tabular}
\end{table}
On the other hand, we know from Lemma \ref{lem23} that $S_{\Gamma_0(14)} := \{\infty, 0, 1/2, 1/7\}$ is the set of all inequivalent cusps of $\Gamma_0(14)$. Moreover, we find that the following triples of cusps of $\Gamma_1(14)$ are equivalent under $\Gamma_0(14)$: $(\infty, 3/14, 5/14), (0,1/3,1/5), (1/2,1/6,1/10)$, and $(1/7,3/7, 5/7)$. Thus, from Table \ref{tab:tbl31} and \cite[Lemma 2]{radu}, we read the orders of vanishing of $gz$ at each element of $S_{\Gamma_0(14)}$ as follows, with $\infty$ now viewed as a cusp of $\Gamma_0(14)$: 
\begin{align*}
	\mathrm{ord}_{\infty}\,gz &= \mathrm{ord}_{3/14}\,gz = \mathrm{ord}_{5/14}\,gz=-4,\\
	\mathrm{ord}_{1/2}\,gz &= \mathrm{ord}_{1/6}\,gz = \mathrm{ord}_{1/10}\,gz \geq 2,\\
	\mathrm{ord}_{0}\,gz &= \mathrm{ord}_{1/3}\,gz = \mathrm{ord}_{1/5}\,gz \geq 0,\\
	\mathrm{ord}_{1/7}\,gz &= \mathrm{ord}_{3/7}\,gz = \mathrm{ord}_{5/7}\,gz \geq 0.
\end{align*}
We deduce that $gz\in\mathcal{M}^\infty(14)$ with unique pole of order $-4$ at $\infty$. We infer from Lemma \ref{lem25} that there is a polynomial
\begin{align*}
	F(X,Y) = X^3-Y^4 +\sum_{\substack{a, b\geq 0\\ 3a+4b\leq 12}}C_{a,b}X^aY^b\in \mathbb{C}[X,Y]
\end{align*}  
such that $F(gz,g^2)=0$. Plugging in the $q$-expansions 
\begin{align}
	z &= q^{-5/2}+2q^{-3/2}+4q^{-1/2}+4q^{1/2}+6q^{3/2}+8q^{5/2}+O(q^3),\nonumber\\
	g &= q^{-3/2}+ 2q^{-1/2}+q^{1/2}+2q^{3/2}+2q^{5/2}+O(q^3)\label{eq32}
\end{align}
yields the polynomial
\begin{align*}
	F(X,Y) = X^3 +4X^2 Y  - 3XY^2 - Y^2(Y^2+4Y+49).
\end{align*}
Since
\begin{align*}
0 = F(gz,g^2) = g^3(z^3+4gz^2-3g^2z-g^5-4g^3-49),
\end{align*}
we arrive at the equation 
\begin{align}\label{eq33}
z^3+4gz^2-3g^2z-g^5-4g^3-49=0,
\end{align}
which is equivalent to (\ref{eq18}).
\end{proof}

To prove Theorem \ref{thm12}, we need the following result on a modular function on $\Gamma_0(14)$ involving 
\begin{align*}
	h(\tau) := \dfrac{\Pi_{q^7}^2}{\Pi_{q^{14}}^2} = \dfrac{\eta^{12}(14\tau)}{\eta^4(7\tau)\eta^8(28\tau)}.
\end{align*}

\begin{lemma}\label{lem31}
The function 
\begin{align*}
H:=H(\tau):= g\left(h(\tau)+\dfrac{16}{h(\tau)}\right)
\end{align*}
lies on $\mathcal{M}^\infty(14)$ with $\mathrm{ord}_\infty\,H = -5$.  
\end{lemma}

\begin{proof}
Observe that
\begin{align*}
h_1(\tau) &:= gh(\tau) = \dfrac{\eta^4(2\tau)\eta^8(14\tau)}{\eta^2(\tau)\eta^2(7\tau)\eta^8(28\tau)},\\
h_2(\tau) &:= \dfrac{h(\tau)}{g} = \dfrac{\eta^2(\tau)\eta^{16}(14\tau)}{\eta^4(2\tau)\eta^6(7\tau)\eta^8(28\tau)}
\end{align*}
are modular on $\Gamma_0(28)$ by Lemma \ref{lem26}. Since $\Gamma_0(14)$ is generated by $\Gamma_0(28)$ and $\alpha:=[\begin{smallmatrix}
	1 & 0\\14 & 1
\end{smallmatrix}]$, we know that $H(\tau) = h_1(\tau)+h_1(\alpha\tau)$ is modular on $\Gamma_0(14)$. We claim that 
\begin{align}\label{eq34}
	h_1(\alpha \tau)=\dfrac{16}{h_2(\tau)}.
\end{align}

Consider the function $J(\tau) := h_1(\alpha\tau)h_2(\tau)$. As $\Gamma_0(28)$ is a normal subgroup of $\Gamma_0(14)$, we have that $J(\tau)$ is modular on $\Gamma_0(28)$. We wish to compute the orders of vanishing of $J(\tau)$ at the cusps of $\Gamma_0(28)$ given by 
\begin{align*}
	S_{\Gamma_0(28)} := \{\infty, 0, 1/2, 1/4, 1/7, 1/14\}
\end{align*}
according to Lemma \ref{lem23}. We look for cusps $r_1, r_2\in S_{\Gamma_0(28)}$ such that $\alpha(r_1)$ is equivalent to $r_2$ over $\Gamma_0(28)$, which will be denoted by $\alpha(r_1)\sim r_2$. In this case, by \cite[Lemma 2]{radu}, we have 
\begin{align}\label{eq35}
	\mbox{ord}_{r_1}\,h_1(\alpha\tau) = \mbox{ord}_{\alpha r_1}\,h_1(\tau) = \mbox{ord}_{r_2}\,h_1(\tau).
\end{align}
From Lemma \ref{lem23}, we get 
\begin{align*}
	&\alpha(0)\sim 0, \qquad \alpha(1/2)\sim 1/4, \qquad \alpha(1/4)\sim 1/2,\\
	&\alpha(1/7)\sim 1/7, \qquad \alpha(1/14)\sim \infty, \qquad \alpha(\infty)\sim 1/14.
\end{align*}

We then compute the required orders using Lemma \ref{lem27}, the above equivalences, and (\ref{eq35}), as shown in Table \ref{tab:tbl32} below. We deduce from this table that $J(\tau)$ is holomorphic at every cusp of $\Gamma_0(28)$, so $J(\tau)$ is constant by \cite[Lemma 5]{radu}. 

\begin{table}[h]
	\caption{The order of vanishing of $h_1(\tau), h_2(\tau)$, and $J(\tau)$ at the cusps of $\Gamma_0(28)$}\label{tab:tbl32}
	\begin{tabular}{@{}ccccccc@{}}
		\toprule
		cusp $r$ & $0$  & $1/2$  & $1/4$ &$1/7$ & $1/14$ & $\infty$  \\
		\midrule
		$\mbox{ord}_r\,h_1(\alpha\tau)$ & $0$ & $1$ & $2$ & $0$ & $-5$ & $2$  \\
		\midrule
		$\mbox{ord}_r\,h_2(\tau)$ & $0$ & $-1$ & $-2$ & $0$& $5$ & $-2$  \\
		\midrule
		$\mbox{ord}_r\,J(\tau)$ & $0$ & $0$ & $0$ & $0$ & $0$ & $0$  \\
		\bottomrule
	\end{tabular}
\end{table}

To find the required constant, we only need to evaluate $J(0)$. In view of 
\begin{align*}
	h_1(\tau)h_2(\tau) = \dfrac{\eta^{24}(14\tau)}{\eta^8(7\tau)\eta^{16}(28\tau)}
\end{align*}
and the transformation formula \cite[Lemma 1]{yang1}
\begin{align*}
	\eta\left(-\dfrac{1}{\tau}\right) = \sqrt{-i\tau}\eta(\tau), 
\end{align*}
we obtain 
\begin{align*}
	J(0) = h_1(0)h_2(0) =\lim_{\tau\rightarrow\infty}\dfrac{\eta^{24}(-14/\tau)}{\eta^8(-7/\tau)\eta^{16}(-28/\tau)}=16\lim_{\tau\rightarrow\infty}\dfrac{(q^{1/14};q^{1/14})_\infty^{24}}{(q^{1/7};q^{1/7})_\infty^8(q^{1/28};q^{1/28})_\infty^{16}}=16,
\end{align*}
proving (\ref{eq34}) as claimed. Thus, $H(\tau)$ is indeed modular on $\Gamma_0(14)$. We next look at Table \ref{tab:tbl32} to get the orders of vanishing of $H(\tau)$ at the cusps of $\Gamma_0(28)$, as shown in Table \ref{tab:tbl33}. We use the inequality
\begin{align*}
	\mathrm{ord}_r\,H(\tau) &\geq \min\{\mathrm{ord}_r\, h_1(\tau),\mathrm{ord}_r\, 1/h_2(\tau)\},
\end{align*}
for $r\in S_{\Gamma_0(28)}$, where equality holds when $\mathrm{ord}_r\, h_1(\tau)\neq \mathrm{ord}_r\, 1/h_2(\tau)$.

\begin{table}[h!]
\caption{The orders of $h_1(\tau), 1/h_2(\tau)$, and $H(\tau)$ at the cusps $r$ of $\Gamma_0(28)$}\label{tab:tbl33}
\begin{tabular}{@{}ccccccc@{}}
	\toprule
	cusp $r$ & $0$  & $1/2$  & $1/4$ &$1/7$ & $1/14$ & $\infty$  \\
	\midrule
	$\mbox{ord}_r\,h_1(\tau)$ & $0$ & $2$ & $1$ & $0$ & $2$ & $-5$  \\
	\midrule
	$\mbox{ord}_r\,1/h_2(\tau)$ & $0$ & $1$ & $2$ & $0$& $-5$ & $2$  \\
	\midrule
	$\mbox{ord}_r\,H(\tau)$ & $\geq 0$ & $1$ & $1$ & $\geq 0$ & $-5$ & $-5$  \\
	\bottomrule
\end{tabular}
\end{table}
Recall that $H(\tau)$ is modular on $\Gamma_0(14)$. In view of Lemma \ref{lem23}, we find the following pairs of cusps of $\Gamma_0(28)$ that are equivalent under $\Gamma_0(14)$, namely $(\infty, 1/14)$ and $(1/2, 1/4)$. We now obtain the orders of $H(\tau)$ at each element of $S_{\Gamma_0(14)}$ using Table \ref{tab:tbl33} and \cite[Lemma 2]{radu}, with $\infty$ now considered as a cusp of $\Gamma_0(14)$:
\begin{align*}
	\mathrm{ord}_{\infty}\,H(\tau) &= \mathrm{ord}_{1/28}\,H(\tau) = -5,\\
	\mathrm{ord}_{1/2}\,H(\tau) &= \mathrm{ord}_{1/4}\,H(\tau) =1,\\
	\mathrm{ord}_{0}\,H(\tau) &\geq 0,\\
	\mathrm{ord}_{1/7}\,H(\tau) &\geq 0.
\end{align*}
Hence, $H(\tau)$ has a unique pole of order $-5$ at $\infty$ and holomorphic elsewhere, so that $H(\tau)\in\mathcal{M}^\infty(14)$ as desired.
\end{proof}

We can now settle Theorem \ref{thm12} using Lemma \ref{lem31}.

\begin{proof}[Proof of Theorem \ref{thm12}]
Substituting $a=q^2, q^4, q^6$ in (\ref{eq31}) and adding the resulting identities gives
\begin{align}
\sum_{n\geq 1}\left[\dfrac{q^{2n}}{(1-q^{2n})^2}-\dfrac{q^{14n}}{(1-q^{14n})^2}\right]-6\sum_{n\geq 1}\dfrac{q^{14n-7}}{(1-q^{14n-7})^2} &= \Pi_{q^7}^2\left(\dfrac{1}{g_1(\tau)}+\dfrac{1}{g_2(\tau)}+\dfrac{1}{g_3(\tau)}\right)\nonumber\\
&=\dfrac{\Pi_{q^7}^2 f_1}{g}\label{eq36}.
\end{align}
Applying the identity 
\begin{align*}
	\sum_{n\geq 1}\left[\dfrac{q^n}{(1-q^n)^2}-\dfrac{q^{2n}}{(1-q^{2n})^2}\right]=\sum_{n\geq 1}\dfrac{q^{2n-1}}{(1-q^{2n-1})^2}
\end{align*}
on (\ref{eq36}) yields
\begin{align}
	g\left(\dfrac{1}{\Pi_{q^7}^2}\sum_{n\geq 1}\left[\dfrac{q^n}{(1-q^n)^2}-7\dfrac{q^{7n}}{(1-q^{7n})^2}\right]-z\right)&-f_1\nonumber\\
	&=\dfrac{6g}{\Pi_{q^7}^2}\sum_{n\geq 1}\left[\dfrac{q^{14n-7}}{(1-q^{14n-7})^2}-\dfrac{q^{14n}}{(1-q^{14n})^2}\right].\label{eq37}
\end{align}
We replace $q$ with $q^7$ in identity (\ref{eq11}), put in on the right hand side of (\ref{eq37}), and then multiply both sides by $4$, arriving at
\begin{align*}
	g(w-4z)-4f_1=H,
\end{align*}
which is equivalent to
\begin{align}\label{eq38}
	gz(f-4) = H+4f_1.
\end{align}
Recall from the proof of Theorem \ref{thm11} that $f_1$ is modular on $\Gamma_0(14)$. We infer from Table \ref{tab:tbl31} and the paragraph immediately after that that the orders of vanishing of $f_1$ at each element of $S_{\Gamma_0(14)}$ are as follows: 
\begin{align*}
	\mathrm{ord}_{\infty}\,f_1 &= \mathrm{ord}_{3/14}\,f_1 = \mathrm{ord}_{5/14}\,f_1=-3,\\
	\mathrm{ord}_{1/2}\,f_1 &= \mathrm{ord}_{1/6}\,f_1 = \mathrm{ord}_{1/10}\,f_1 \geq 1,\\
	\mathrm{ord}_{0}\,f_1 &= \mathrm{ord}_{1/3}\,f_1 = \mathrm{ord}_{1/5}\,f_1 \geq 0,\\
	\mathrm{ord}_{1/7}\,f_1 &= \mathrm{ord}_{3/7}\,f_1 = \mathrm{ord}_{5/7}\,f_1 \geq 0.
\end{align*}
Thus, we see that $f_1\in \mathcal{M}^\infty(14)$ with unique pole of order $-3$ at $\infty$. Furthermore, $t:= H+4f_1$ also lies on $\mathcal{M}^\infty(14)$ with $\mathrm{ord}_\infty\,t=-5$. Lemma \ref{lem25} asserts that there is a polynomial
\begin{align*}
	G_1(X,Y) = X^3-Y^5 +\sum_{\substack{a, b\geq 0\\ 3a+5b\leq 15}}C_{a,b}X^aY^b\in \mathbb{C}[X,Y]
\end{align*}
such that $G_1(gz(f-4),g^2) = G_1(t,g^2)=0$ from (\ref{eq38}). Substituting the $q$-expansions (\ref{eq32}) and 
\begin{align*}
	t = q^{-5}+2q^{-4}+5q^{-3}+10q^{-2}+18q^{-1}+32+O(q),
\end{align*}
we find that
\begin{align*}
	G_1(X,Y) = X^3 - 33YX^2+2Y(7Y^2+46Y+343)X - Y(Y^2+2Y+49)^2.
\end{align*}
We recast (\ref{eq33}) as $G_2(z,g) = 0$, where 
\begin{align*}
	G_2(X,Y):= X^3+4X^2Y-3XY^2 - Y(Y^4+4Y^2+49).
\end{align*}
We now eliminate $z$ by applying the resultant of $G_1(gz(f-4),g^2)$ and $G_2(z,g)$ viewed as polynomials in $z$. We deduce that
\begin{align}\label{eq39}
0 &= \mathrm{Res}_z\left(G_1(gz(f-4),g^2), G_2(z,g)\right) = P(f,g)Q(f,g), 
\end{align}
where
\begin{align*}
P(X,Y) &:= Y^2(Y^4+4Y^2+49)X^3 - Y^2(2Y^4+5Y^2+98)X^2 - 2Y^2(5Y^4+22Y^2+245)X\\
&- (Y^2+4Y+7)^2(Y^2-4Y+7)^2,\\
Q(X,Y) &:= Y^4(Y^4+4Y^2+49)^2X^6 - 2Y^4(Y^4+4Y^2+49)(17Y^4+20Y^2+833)X^5\\
&+ 2Y^4(315Y^8-886Y^6+22304Y^4-43414Y^2+756315)X^4\\
& -2Y^2(Y^{12}+3097Y^{10}-26139Y^8+150074Y^6-1280811Y^4+7435897Y^2+117649)X^3\\
&+ Y^2(34Y^{12}+33705Y^{10}-359550Y^8+1335446Y^6-17617950Y^4+80925705Y^2\\
&+4000066)X^2-2Y^2(49Y^{12}+47722Y^{10}-505249Y^8+1923820Y^6-24757201Y^4\\
&+114580522Y^2+5764801)X + Y^{16}-16Y^{14}+108572Y^{12}-1032688Y^{10}+4938886Y^8\\
&-50601712Y^6+260681372Y^4-1882382Y^2+5764801.
\end{align*}
Substituting the $q$-expansions (\ref{eq32}) and 
\begin{align*}
	f = \dfrac{1}{q}+2+4q-4q^2+6q^3-8q^4+O(q^5),
\end{align*}
we verify that $Q(f,g)\neq 0$. Hence, we conclude from (\ref{eq39}) that $P(f,g)=0$, which yields (\ref{eq19}) as desired.
\end{proof}

%
%
%
%
%

\section*{Acknowledgments}
The author would like to thank Dr. M. V. Yathirajsharma for sending him a copy of the paper \cite{yat}. 

\bibliographystyle{amsplain}
\bibliography{gosper14}
\end{document}